\documentclass[twoside]{article}
\usepackage{amssymb,a4}
\usepackage{amsthm}
\def\cen{\centerline}

\def\ad{\hbox{ad}}

\def\a{\alpha}
\def\b{\beta}

\def\ll{\lambda}

\def\gl{\frak{gl}}
\newfont{\df}{eufm10}
\def\vep{\varepsilon}

\def\ll{\lambda}
\def\der{\hbox{Der}\,}
\def\aut{\hbox{Aut}\,}
\def\fsl{\frak{sl}\,}
\def\fst{\frak{st}\,}
\def\stl{\frak{stl}\,}

\def\vp{\varphi}

\def\r{\gamma}

\def\ot{\otimes}

\def\de{\delta}
\def\De{\Delta}

\def\ker{\hbox{Ker}\,}

\def\ad{\hbox{\rm ad}\,}
\def\Ad{\hbox{\rm Ad}\,}

\def\inn{\hbox{\rm Inn}\,}

\def\ot{\otimes}

\def\mg{{\bf \frak g}}
\def\dmg{\dot{\mg}}

\hoffset
\voffset
\title{{Leibniz Algebras Graded by Finite Root Systems}
\footnotetext{E-mail: *liudong@hutc.zj.cn,
**nhhu@euler.math.ecnu.edu.cn}}
\author{{\bf Dong Liu*} \\ Department of Mathematics, Huzhou Teachers College\\ Zhejiang Huzhou, 313000, P.R. China\\ \\{\bf Naihong Hu**}\\Department of Mathematics, East China Normal University\\ Shanghai, 200062, P.R. China }
\date{ }

\begin{document}
\maketitle
\def\abstractname{ABSTRACT}
\begin{abstract}
There are several researches on Lie algebras and Lie superalgebras
graded by finite root systems. In this paper, we study Leibniz
algebras graded by finite root systems and obtain some important
results in simply-laced cases.
\smallskip


\end{abstract}

\newtheorem{theo}{Theorem}[section]
\newtheorem{defi}[theo]{Definition}
\newtheorem{lemm}[theo]{Lemma}
\newtheorem{coro}[theo]{Corollary}
\newtheorem{prop}[theo]{Proposition}

\section{Introduction}

In \cite{Lo2}, J.-L. Loday introduced a non-antisymmetric version of
Lie algebras, whose bracket satisfies the Leibniz relation (see
(2.5)), therefore  called {\it Leibniz algebra}. The Leibniz
relation, combined with antisymmetry, is a variation of the Jacobi
identity, hence Lie algebras are anti-symmetric Leibniz algebras. In
\cite{Lo4}, Loday also introduced an 'associative' version of
Leibniz algebras, called {\it associative dialgebras}, equipped with
two binary operations, $\vdash$ and $\dashv$, which satisfy the five
relations (see the axiom (Ass) in section 2). These identities are
all variations of the associative law, so associative algebras are
dialgebras for which the two products coincide. The peculiar point
is that the bracket $[a, b]=:a\dashv b-b\vdash a$ defines a Leibniz
algebra which is not antisymmetric, unless the left and right
products coincide. Hence dialgebras yield a commutative diagram of
categories and functors
\begin{eqnarray*}
{\bf Dias}&\stackrel{-}{\to}& {\bf Leib}\\
\downarrow&&\downarrow\\
{\bf Assoc}&\stackrel{-}{\to}& {\bf Lie}
\end{eqnarray*}

Steinberg Lie algebras come from Steinberg groups, which are closely
connected with K-theory, and play a key role in the study of Lie
algebras graded by finite root systems of type $A$. By definition,
the {\it Steinberg Lie algebra} $\fst(n, A)$ over a $K$-algebra $A$
is a Lie algebra generated by symbols $v_{ij}(a)$, $1\le i\ne j\le
n$, $a\in A$, subject to the relations
\begin{eqnarray}
&&v_{ij}(k_1a+k_2b)=k_1v_{ij}(a)+k_2v_{ij}(b),\  \hbox{ for } \ a, b\in D, \ k_1,
k_2\in K;\\
&&[v_{ij}(a), v_{kl}(b)]=0\  \hbox{ if }\ i\ne l\ \hbox{ and } \  j\ne k;\\
&&[v_{ij}(a), v_{kl}(b)]=v_{il}(ab)\  \hbox{ if } \ i\ne l\
\hbox{and } \  j= k.
\end{eqnarray}
It is clear that the relation (3) makes sense only if $n\ge3$.

 From \cite{Fa} we see that the map $\eta: a\to v_{ij}(a)$ is one-to-one if and
 only if $A$ is an associative algebra for $n\ge 4$ and  $A$ is an alternative
 algebra for $n=3$.

In 1992, S. Berman and R.V. Moody (\cite{BM}) studied Lie algebras
graded by finite root systems of $A_l\,(l\ge 2)$, $D_l\, (l\ge 4)$,
$E_l\, (l=6, 7, 8)$ and obtained the structure of a Lie algebra over
$K$ graded by the root system $\Delta$ of type $X_l\, ( l\ge2)\, (
X_l=A_l, D_l, E_l)$.

The universal central extensions of Lie algebras graded by finite
root systems were studied in several papers (\cite{B}, \cite{Gar},
\cite{KL}, \cite{Gao1}, \cite{ABG1}, etc.).

 In this paper we shall consider Leibniz algebras graded by finite root systems of
 types $A, D$ and $E$. We also prove that
\begin{theo}{\bf(Recognition Theorem).}
 Let $L$ be a Leibniz algebra over $K$ graded by the root system $\Delta$ of
 type $X_l(l\ge2)\, (X_l=A_l, D_l, E_l)$.

(1) If $X_l=A_l, l\ge 3$, then there exists a unital associative
$K$-dialgebra $R$ such that $L$ is centrally isogenous with
$\fsl(l+1, R)$;

(2)  If $X_l=A_l, l=2$, then there exists a unital alternative
$K$-dialgebra $R$ such that $L$ is centrally isogenous with
$\stl(l+1, R)$, where $\stl(n, R)$ is defined in Section 2.4;

(3)  If $X_l=D_l\, (l\ge 4), E_l\, (l=6, 7, 8)$, then there exists a
unital associative commutative $K$-dialgebra $R$ such that $L$ is
centrally isogenous with $\dmg\ot R$.

\end{theo}
{\bf Remark.} Two perfect Lie algebras $L_1$ and $L_2$ are called
{\it centrally isogenous} if they have the same universal central
extension (up to isomorphism).

The paper is organized as follows. In Section 2, we recall some
notions of Leibniz algebras and dialgebras. In Section 3, we give
the definition of Leibniz algebras graded by finite root systems. In
Sections 4 and 5, we mainly prove the Recognition Theorem (Theorem
1.1). Throughout this paper, $K$ denotes a field of characteristic
0, $R$ a unital dialgebra over $K$.

\section{Dialgebras and Leibniz algebras}

We recall the notions of associative dialgebras, alternative
dialgebras, Leibniz algebras and their (co)homology as defined in
\cite{Lo1}---\cite{Lo4} and \cite{L}.

\subsection{Dialgebras.}

\begin{defi}\cite{Lo4}  A {\it dialgebra} $D$ over $K$ is a $K$-vector space
$D$ with two operations $\dashv, \vdash:D\ot D\to D$, called left
and right products.
\end{defi}

A dialgebra is called unital if it is given a specified bar-unit: an
element $1\in D$ which is a unit for the left and right products
only on the bar-side, that is, $1\vdash a=a=a\dashv 1$, for any
$a\in D$. A morphism of dialgebras is a $K$-linear map $f:D\to D'$
which preserves the products, i.e. $f(a\star b)=f(a)\star f(b)$,
where $\star$ denotes either the product $\dashv$ or the product
$\vdash$.

\begin{defi} \cite{Lo4} A dialgebra $D$ over $K$ is called {\it associative}
if the two operators $\dashv$ and $\vdash$ satisfy the following
five axioms:
$$\cases{a\dashv(b\dashv c)=(a\dashv b)\dashv c=a\dashv(b\vdash c),\cr
        (a\vdash b)\dashv c=a\vdash(b\dashv c),\cr
         (a\vdash b)\vdash c=a\vdash (b\vdash c)=(a\dashv b)\vdash c.  }\leqno(Ass)$$
\end{defi}

Denote by {\bf Dias, Assoc} the categories of associative dialgebras
and associative algebras over $K$ respectively. Then the category
{\bf Assoc} is a full subcategory of {\bf Dias}.

Obviously, an associative dialgebra is an associative algebra if and
only if $a\dashv b=a\vdash b=ab$.

The concept of alternative dialgebras was introduced in \cite{L} for
the study of the Steinberg Leibniz algebras.

\begin{defi}\cite {L} A dialgebra $D$ over $K$ is called {\it alternative}
if the two operators $\dashv$ and $\vdash$ satisfying the following five axioms:
$$\cases{J_{\dashv}(a, b, c)=-J_{\vdash}(c, b, a), \quad J_{\dashv}(a, b, c)=J_{\vdash}(b, c, a),\cr
        J_{\times}(a, b, c)=-J_{\vdash}(a, c, b),\cr
         (a\vdash b)\vdash c=(a\dashv b)\vdash c,\quad  a\dashv (b\vdash c)
         =a\dashv (b\dashv c), }\leqno(Alt)$$
where $J_{\dashv}(a, b, c)=(a\dashv b)\dashv c-a\dashv(b\dashv c),\  J_{\vdash}(a, b, c)
=(a\vdash b)\vdash c-a\vdash(b\vdash c)$ and $J_{\times}(a, b, c)
=(a\vdash b)\dashv c-a\vdash(b\dashv c)$.
\end{defi}

Obviously, an alternative dialgebra is an alternative algebra if
$a\dashv b=a\vdash b=ab$. Moreover, the following formulae are clear
for an alternative dialgebras according to the definition.
$$J_{\dashv}(a, b, c)=-J_{\dashv}(a, c, b),\eqno(2.1)$$
$$J_{\vdash}(a, b, c)=-J_{\vdash}(b, a, c),\eqno(2.2)$$
$$J_{\times}(a, b, c)=-J_{\times}(c, b, a).\eqno(2.3)$$
So we also have $$J_{\dashv}(a, b, b)=0, \  J_{\vdash}(a, a, b)=0, \
J_{\times}(a, b, a)=0.\eqno(2.4)$$ {\bf Examples.}

1. Obviously, an associative (alternative) dialgebra is an
associative (alternative) algebra if and only if $a\dashv b=a\vdash
b=ab$.

2. {\it Differential associative (alternative) dialgebra.} Let $(A, d)$ be
a differential associative (alternative) algebra. So by hypothesis,
$d(ab)=(da)b+adb$ and $d^2=0$. Define left and right products on
$A$ by the formulas
$$x\dashv y=xdy, \quad x\vdash y=(dx)y.$$
Then $A$ equipped with these two products is an associative (alternative) dialgebra.

3. Tensor product. Let $D$ and $D'$ be two associative dialgebras,
then $D\ot D'$ with multiplication $(a\ot a')\star (b\ot b')=(a\star
b)\ot (a'\star b')$, $\star=\dashv, \vdash$, is also an associative
dialgebra. Especially, if $D$ is a unital associative dialgebra,
then $M_n(D)=M_n(K)\ot D$ is also a unital associative dialgebra.

4. Let $D$ be an associative (alternative) algebra.  On the module of $n$-space $D=A^{\otimes n}$ one puts
$$(x\dashv y)_i=x_i(\sum_{j=1}^ny_j), \  i=1, \cdots, n\quad \hbox {and}$$
$$(x\vdash y )_i=(\sum_{j=1}^nx_j)y_i, \  i=1, \cdots, n. $$
Then $(D, \dashv, \vdash)$ is an associative (alternative) dialgebra. For $n=1$, this is example 1.

\subsection{Leibniz algebras.} A {\it Leibniz algebra} \cite {Lo2} $L$
is a vector space over a field $K$ equipped with a $K$-bilinear
map
$$[-,-]: L\times L\to L$$
satisfying the Leibniz identity
$$[x, [y, z]]= [[x, y], z]-[[x, z], y], \quad \forall \;x, \,y, \,z\in L.\eqno(2.5)$$

Obviously, a Lie algebra is a Leibniz algebra. A Leibniz algebra is a Lie algebra if
and only if
$[x, x]=0$ for all $x\in L$.

Suppose that $L$ is a Leibniz algebra over $K$. For any $z\in L$,
we define $\ad z\in \hbox{End}_kL$ by
$$\ad z(x)=-[x, z], \quad\forall x\in L.$$
It follows (2.5) that
$$\ad z([x, y])=[\ad z(x), y]+[x, \ad z(y)]$$
for all $x, y\in L$. This says that $\ad z$ is a derivation of
$L$. We also call it an inner derivation of $L$.

Similarly, we also have the definition of general derivation of a
Leibniz algebra and we denote by $\inn(L)$, $\der(L)$ the set of all
inner derivations, derivations of $L$ respectively. They are also
Leibniz algebras.

Let $L$ be a Leibniz algebra over $K$. Consider the boundary map: $\de_n:L^{\ot n}\to L^{\ot (n-1)}$ defined by
$$\de_n(x_1\ot\cdots\ot x_n)=\sum_{1\le i<j\le n}(-1)^{j+1}x_1\ot \cdots \ot x_{i-1}\ot[x_i, x_j]\ot x_{i+1}\ot\cdots\ot\hat x_j\ot\cdots\ot x_n,$$
where $\hat x_j$ indicates that the term $x_j$ is omitted. One can
show that $\de^2=0$ (see \cite{LP1}) and the complex $(L^{\ot
n},\, \de)\ (L^0=K,\, \de_1=0)$ gives the Leibniz homology
$HL_*(L)$ of the Leibniz algebra $L$.

Let $L$ be a Leibniz algebra over $K$. It is called perfect if
$[L, L]=L$. A central extension of $L$ is a pair $(\hat L, \pi)$
where $\hat L$ is a Leibniz algebra and $\pi: \hat L\to L$ is a
surjective homomorphism such that $\ker\pi$ lies in the center of
$\hat L$ and the exact sequence $0\to \ker\pi\to \hat L\to L\to 0$
splits as $K$-module. The pair $(\hat L, \pi)$ is a universal
central extension of $L$ if for every central extension $(\tilde
L, \tau)$ of $L$ there is a unique homomorphism $\psi:\hat L\to
\tilde L$ for which $\tau\circ\psi=\pi$. So the universal central
extension is unique, up to isomorphism. A Leibniz algebra $L$ has
a universal central extension if and only if $L$ is perfect. If
$(\hat L, \pi)$ is the universal central extension of $L$, then
$$HL_2(L)\cong \ker\pi.\eqno(2.6)$$

\noindent{\bf Remark.} In \cite{LP1}, \cite{Gao1}, \cite{Gao2} and
\cite{LH1}, the universal central extensions of  many infinite
dimensional Lie algebras in the category of Leibniz algebras are
determined.

We also denote by {\bf Leib} and {\bf Lie} the categories of Leibniz
algebras and Lie algebras over $K$ respectively.

For any associative dialgebra $D$, define $$[x, y]=x\dashv y-y\vdash x,$$
then $D$ equipped with this bracket is a Leibniz algebra. We denote it by $D_L$.
 The canonical map $D\mapsto D_{L}$ induces a functor $(-):$ {\bf Dias}$\to${\bf Leib}.

For a Leibniz algebra $L$, let $L_{Lie}$ be the quotient of $L$ by
the ideal generated by the elements $[x, y]+[y, x]$, for all $x,
y\in L$. It is clear that $L_{Lie}$ is a Lie algebra. The canonical
projection $L\to L_{Lie}$ is universal among the maps from $L$ to
Lie algebras. In other words, the functor $(-)_{Lie}:$ {\bf
Leib}$\to${\bf Lie} is left adjoint to $inc:$ {\bf Lie}$\to${\bf
Leib}.

Moreover, we have the following commutative diagram of categories
and functors
\begin{eqnarray*}
{\bf Dias}&\stackrel{-}{\to}& {\bf Leib}\\
\downarrow&&\downarrow\\
{\bf Assoc}&\stackrel{-}{\to}& {\bf Lie}
\end{eqnarray*}

As in the Lie algebra case, the universal enveloping associative
dialgebra of a Leibniz algebra $L$ is defined as
$$Ud(L):=(T(L)\ot L\ot T(L))/\{[x, y]-x\dashv y+y\vdash x|x, y\in L\},$$
where elements $x,y$ of $L$ are regarded as elements in $K\otimes
L\otimes K$.

\begin{prop}$\cite{Lo4}$

 The functor $Ud:{\bf Leib}\to {\bf Dias}$ is left adjoint to the functor $-:{\bf Dias}\to {\bf Leib}$.
\end{prop}

Let $L$ be a Leibniz algebra, then $M$ is said to be a right $L$-module if $M$ is a $K$-vector space  equipped with the action of $L$:
$$[-,-]: M\times L\to M$$
satisfying

$$[m,[x,y]]=[[m,x],y]-[[m,y],x],\ \forall x,\, y\in L, m\in M.$$

So for a Lie algebra $\mg$, any right $\mg$-module in the Leibniz algebra case is just the right $\mg$-module in the Lie algebra case.

\subsection{Lie algebras graded by finite root system.}

First we introduce the Steinberg Lie algebra. Steinberg Lie algebras
come from Steinberg groups, which are closely connected with
K-theory.

If $R$ is a (nonassociative) ring with 1, and $n\ge3$, the {\it
Steinberg group} (see \cite{Fa}) St$_n(R)$ is the group generated by
the symbols $x_{ij}(a), 1\le i\ne j\le n, a\in R$, subject to the
relations
\begin{eqnarray*}
&&x_{ij}(a+b)=x_{ij}(a)x_{ij}(b),\  \hbox{ for } \ a, b\in R, \ k_1, k_2\in K;\\
&&[x_{ij}(a), x_{kl}(b)]=1\  \hbox{ if }\ i\ne l\ \hbox{ and } \  j\ne k;\\
&&[x_{ij}(a), x_{kl}(b)]=x_{il}(ab)\  \hbox{ if } \ i\ne l\
\hbox{and } \  j= k,
\end{eqnarray*} where $[x, y]=xyx^{-1}y^{-1}$ is the group commutator.

By definition, the {\it Steinberg Lie algebra} (see \cite{Fa})
$\fst(n, A)$ $(n\ge 3)$ over a $K$-algebra $A$ is a Lie algebra
generated by symbols $u_{ij}(a)$, $1\le i\ne j\le n$, $a\in A$,
subject to the relations
\begin{eqnarray*}
&&u_{ij}(k_1a+k_2b)=k_1u_{ij}(a)+k_2u_{ij}(b),\  \hbox{ for } \ a, b\in A, \ k_1, k_2\in K;\\
&&[u_{ij}(a), u_{kl}(b)]=0\  \hbox{ if }\ i\ne l\ \hbox{ and } \  j\ne k;\\
&&[u_{ij}(a), u_{kl}(b)]=u_{il}(ab)\  \hbox{ if } \ i\ne l\
\hbox{and } \  j= k.
\end{eqnarray*}

Define $i_n(A)=\{a\in A\mid u_{ij}(a)=0 \}$. Clearly, it is an ideal
of $A$ and does not depend on the choice of $i\ne j$.

 From \cite{Fa}, we see that $i_n(A)=0$ if and only if $A$ is an associative algebra
 for $n\ge 4$ and  $A$ is an alternative algebra for $n=3$.

By definition, a $K$-algebra $A$ is called {\it alternative} if $(a, b,
c)=-(c, b, a)=(b, c, a)$, where $(a, b, c)=(ab)c-a(bc)$.

A Lie algebra $L$ over a field $K$ of characteristic 0 is {\it
graded by the (reduced) root system} (see \cite{BM}) $\Delta$ or
is {\it $\Delta$-graded} if

(1) L contains as a subalgebra a finite-dimensional simple Lie
algebra $\dmg=H\oplus\bigoplus_{\a\in\Delta}\dmg_{\a}$ whose root
system is $\Delta$ relative to a split Cartan subalgebra $H=\dmg_0$;

(2) $L=\bigoplus_{\a\in\Delta\cup\{0\}}L_{\a}$, where
$L_{\a}=\{x\in L\mid [h, x]=\a(h)x, \forall h\in H\}$ for $\a\in \Delta\cup\{0\}$; and

(3) $L_0=\sum_{\a\in\Delta}[L_\a, L_{-\a}]$.

\begin{theo} \cite{BM}
 Let $L$ be a Lie algebra over $K$ graded by the root system $\Delta$ of
 type $X_l\, ( l\ge2)\, ( X_l=A_l, D_l, E_l)$ .

(1) If $X_l=A_l,\, l\ge 3$, then there exists a unital associative
$K$-algebra $A$ such that $L$ is centrally isogenous with
$\fsl(l+1, A)$.

(2)  If $X_l=A_l,\, l=2$, then there exists a unital alternative
$K$-algebra $A$ such that $L$ is centrally isogenous with
$\frak{st}(l+1, A)$.

(3)  If $X_l=D_l\,(l\ge 4), \, E_l\,(l=6, 7, 8)$, then there
exists a unital associative commutative $K$-algebra $A$ such that
$L$ is centrally isogenous with $\dmg\ot A$.
\end{theo}

\subsection{Steinberg Leibniz algebras }

The matrix Leibniz algebra
$\gl(n, D)$ is generated by all $n\times n$ matrices with coefficients from a unital
associative dialgebra $D$, and $n\ge 3$ with the bracket
$$[E_{ij}(a), E_{kl}(b)]=\de_{jk}E_{il}(a\dashv b)-\de_{il}E_{kj}(b\vdash a),$$ for all $a, b\in D$, where $E_{ij}(a)$ is the $n\times n$ matrix with coefficient $a$ on $(i, j)$-th position and 0 in all others.

Clearly, $\gl(n, D)$ is a Leibniz algebra. If $D$ is an associative algebra, then  $\gl(n, D)$ becomes a Lie algebra.

Now we consider the subalgebra $\fsl(n, D):=[\,\gl(n, D), \gl(n, D)\,]$, which is called the {\it special linear Leibniz algebra} with coefficients in $D$, of $\gl(n, D)$.

By definition, the {\it special linear Leibniz algebra} $\fsl(n, D)$
has generators $E_{ij}(a),\, 1\le i\ne j\le n, a\in D$, which
satisfy the following relations:
\begin{eqnarray*}
&&[E_{ij}(a), E_{kl}(b)]=0\  \hbox{if}\ i\ne l\ \hbox{and}\  j\ne k;\\
&&[E_{ij}(a), E_{kl}(b)]=E_{il}(a\dashv b)\  \hbox{if}\ i\ne l\ \hbox{and}\  j= k;\\
&&[E_{ij}(a), E_{kl}(b)]=-E_{kj}(b\vdash a)\ \hbox{if}\ i= l\
\hbox{and}\  j\ne k.
\end{eqnarray*}

The Steinberg Leibniz algebra was first introduced in \cite{LP1} for
associative algebras and in \cite{L} for associative dialgebras. By
definition, the {\it Steinberg Leibniz algebra} $\stl(n, D)$ is a
Leibniz algebra generated by symbols $v_{ij}(a)$, $1\le i\ne j\le
n$, $a\in D$, subject to the relations
\begin{eqnarray}
&&v_{ij}(k_1a+k_2b)=k_1v_{ij}(a)+k_2v_{ij}(b),\  \hbox{ for } \ a, b\in D, \ k_1, k_2\in K;\\
&&[v_{ij}(a), v_{kl}(b)]=0,\  \hbox{ if }\ i\ne l\ \hbox{ and } \  j\ne k;\\
&&[v_{ij}(a), v_{kl}(b)]=v_{il}(a\dashv b)\  \hbox{ if } \ i\ne l\ \hbox{and } \  j= k;\\
&&[v_{ij}(a), v_{kl}(b)]=-v_{kj}(b\vdash a)\  \hbox{ if } \ i= l\
\hbox{and} \  j\ne k.
\end{eqnarray}
It is clear that the relations (6)---(7) make sense only if
$n\ge3$.

Let $H_{ij}(a, b):=[v_{ij}(a), v_{ji}(b)]$ for $1\le i\ne j\le n, a,
b\in D$, and $H$ the submodule of $\stl(n, D)$ generated by
$H_{ij}(a, b), i\ne j, a, b\in D$. Define $i_n(D)=\{a\in D\mid
v_{ij}(a)=0 \}$. Clearly, it is an ideal of $D$ and does not depend
on the choice of $i\ne j$. The same consideration as in Steinberg
Lie algebras, we have
\begin{prop}\cite{L}
 For a unital dialgebra $D$, $i_n(D)=0$ in $\stl(n, D)$ if and only if $D$ is
 associative for $n\ge 4$ and  $D$ is alternative for $n=3$.
\end{prop}

The Steinberg Leibniz algebra $\stl(n, D)$ with $n\ge3$ is perfect.

The homomorphism $\psi$ of Leibniz algebras
$$\psi: \quad \stl(n, D)\to \fsl(n, D)$$ by the rule
$\psi(v_{ij}(a))=E_{ij}(a)$ is a surjective homomorphism.

\begin{theo}\cite{L}
If $n\ge3$, then $(\stl(n, D), \psi)$ is the universal central
extension of the Leibniz algebra $\fsl(n, D)$ with kernel $HHS_1(D)$
for a unital associative dialgebra $D$, where $HHS_1(D)$ is the
first homology group of chain complex $(CS_*(D), d)$ defined in
\cite{Fa} (or see \cite{L}) .
\end{theo}

\section{ Leibniz algebras graded by finite root systems}

\begin{defi}
A Leibniz algebra $L$ over a field $K$ of characteristic 0 is
graded by the $($reduced$)$ root system $\Delta$ or is
$\Delta$-graded if

$(1)$ L contains as a subalgebra a finite-dimensional simple Lie
algebra $\dmg=H\oplus\bigoplus_{\a\in\Delta}\dmg_{\a}$ whose root
system is $\Delta$ relative to a split Cartan subalgebra $H=\dmg_0$;

$(2)$ $L=\bigoplus_{\a\in\Delta\cup\{0\}}L_{\a}$, where
$L_{\a}=\{x\in L\mid \ad h(x)=-[x, h]=\a(h)x, \forall h\in H\}$
for $\a\in \Delta\cup\{0\}$; and

$(3)$ $L_0=\sum_{\a\in\Delta}[L_\a, L_{-\a}]$.
\end{defi}

{\bf Remarks.}

1. The conditions for being a $\De$-graded Leibniz algebra imply
that $L$ is a direct sum of finite-dimensional irreducible right
$\dmg$-modules whose highest weights are roots, hence either the
highest long root or short root or $0$.

2. If $L$ is $\Delta$-graded, then $L$ is perfect. Indeed,  the
result follows from $L_\a=[L_\a, H]$ for all $\a\in\Delta$ and (3)
as above.

3. The Steinberg Leibniz algebra $\stl(n, D)$ is graded by the
root system of type $A_{n-1}$. Let $D$ be a commutative dialgebra,
then the Leibniz algebra $\dmg\ot D$ and its central extensions
are graded by the root system of type $\dmg$ (see \cite{LL}).

Now we shall prove the Recognition Theorem (Theorem 1.1). So from
now on, we always set $\Delta$ to be the root system of type $A_l\,
(l\ge2)$, $D_l\, (l\ge 4)$ or $E_l \,(l=6, 7, 8)$.

Such as that in \cite{BM}, we also have the following results.
\begin{defi}
An ordered pair $(\b, \r)\in \Delta\times\Delta$ is an $A_2$-pair if
$(\b, \r)=-1$. Thus $(\b, \r)$ is an $A_2$-pair if and only if it is
a base for an $A_2$ subroot system
 of $\Delta$.
 Two $A_2$-pairs $(\b, \r)$ and $(\b', \r')$ are equivalent,
 written  $(\b, \r)\sim(\b', \r')$,
 if there is an element $w$ of the Weyl group $W$ of $\Delta$ such that $\b'=w\b$
 and $\r'=w\r$.
The equivalent class of $(\b, \r)$ is denoted by $[(\b, \r)]$.
\end{defi}

\begin{lemm}\cite{BM}
$(1)$ If $\Delta$ is of type $D$ or $E$, then there is only one
equivalent class of $A_2$-pairs.

$(2)$ If $\Delta$ is of type $A$, then there are exactly two
equivalent classes of $A_2$-pairs and furthermore, if $(\b, \r)$ is
an $A_2$-pair, then
$$(\b, \r)\sim(-\r, -\b), \  (\b, \r)\not\sim(\r, \b).$$

$(3)$ In all cases if $(\b, \r)$ and $(\r, \de)$ are $A_2$-pairs
with $(\b\mid\de)=0$, then
$$(\b, \r)\sim(\r, \de)\sim(\b, \r+\de)\sim(\b+\r, \de).\quad  \rule[-.23ex]{1.0ex}{2.0ex}$$
\end{lemm}

{\bf Remark.} In type $A_l$, for definiteness, we distinguish the
two classes as follows. We choose a base $\Pi=\{\a_1, \cdots,
\a_l\}$ for $\Delta$ once and for all with Coxeter-Dynkin diagram
$\a_1\circ$----$\circ$----$\cdots$----$\circ\a_l$. Then $[(\a_1,
\a_2)]$ is the positive class and an $A_2$-pair $(\a, \b)\in [\a_1,
\a_2]$ is called positive pair, $[(\a_2, \a_1)]$ is the negative
class and an $A_2$-pair $(\a, \b)\in [\a_2, \a_1]$ is called
negative pair. For convenience, any $A_2$-pair in types $D_l$ and
$E_l$ is either positive and negative.

Let $L$ be a Leibniz algebra graded by $\Delta$ and $\dmg\subset L$
the split simple Lie algebra of Definition 3.1. Let $\{e_\a, H_i\mid
\a\in\Delta, i=1, \cdots, l\}$ be a Chevalley basis of $\dmg$ so
that $e_\a\in L_\a$ for any $\a\in\Delta$, and $\{H_i, i=1, \cdots,
l\}$ a basis of $H$. Let $G=\langle \exp te_\a\mid t\in K\rangle$ be
the corresponding simply connected Chevalley group. For each
$\a\in\Delta$, $\{e_\a,\, \a^{\vee}=[e_\a, e_{-\a}],\, e_{-\a}\}$ is
an $\fsl_2$-triplet. Let
$$n_\a(t)=\exp te_\a\exp(-t^{-1}e_{-\a})\exp te_\a$$ and
set $$N=\langle n_\a(t)\mid \a\in\Delta,\, t\in K^{\times}\rangle,$$
where $K^{\times}$ is the set of non-zero elements of $K$.

Let $h_\a(t)=n_\a(t)n_\a(1)^{-1}$ and $T=\langle h_\a(t)\mid
\a\in\Delta,\, t\in K^{\times}\rangle$. Then $$T\lhd N,\, N/T\cong
W.$$

Clearly, from the $\Delta$-grading $\ad e_\a$ and  $\ad e_{-\a}$ are
nilpotent on $L$, so we can define a homomorphism
$$\Ad: G\to \aut(L)\quad \hbox{by}\quad \exp te_\a\to \exp\ad te_\a, \  \a\in \Delta,\, t\in K.$$

Next, recall that if $M$ is any integrable (right) $\dmg$-module
(one on which $e_\a$ acts locally nilpotent for all $\a\in\Delta$)
with weight space decomposition $M=\oplus M_\ll$ relative to $\dmg$
then, letting $r_\a$ denote the reflection $r_\a\ll=\ll-\langle \ll,
\a^{\vee}\rangle\a$, we have
$$n_\a(t)M_\ll=M_{r_\a\ll} \quad \hbox{for all weights}\ \ll,\ \hbox{and} \eqno(3.1)$$
$$\hbox{the action of}\  h_\a(t)\ \hbox{restricted to}\  M_\ll\  \hbox{is a
scalar-multiplication by}\  t^{\langle \ll, \a^{\vee}\rangle}. \eqno(3.2)$$

In particular, the adjoint representation restricted to $\dmg$ acts
on $L$ in this way.

Fix any $\a\in\Delta$. Let $W^\a$ denote the stabilizer of $\a$ in
$W$ and let $N^\a:=W^{\a}T$. Then
$$W^\a=\langle r_\b\mid \b\in\Delta, (\b|\a)=0\rangle, \eqno(3.3)$$
$$N^\a=\langle n_\b(1)\mid \b\in\Delta, (\b|\a)=0\rangle\cdot T.\eqno(3.4)$$

Fix any $\b\in \Delta$ and choose $w\in W$ with $w\a=\b$. Choose any
$n\in N$ with $nT\leftrightarrow w$ in the isomorphism $N/T\cong W$.
Then $\Ad(n)L_\a=L_\b$ by (3.1). Let the restriction of $\Ad(n)$ to
$L_\a$ be denoted by $\tilde\ll_{\b, \a}$. If also $w'\a=\b$ and
$n'T\leftrightarrow w'$, then $w^{-1}w'\in W^\a$ and $n^{-1}n'\in
N^\a$. Thus $n'\in nN^\a=n\langle n_\r(1)\mid \r\in\Delta,\,
(\b|\a)=0\rangle T$. Elements of $T$ acts as scalar multiplications
on $L_\a$. Furthermore,
$$(\r|\a)=0\Rightarrow (\a\pm\r, \a\pm\r)=4\Rightarrow \a\pm\r\not\in\Delta\cup\{0\}$$
$$\Rightarrow \Ad\exp(te_{\pm\r})\ \hbox{acts as the identity on}\ L_\a, \forall\, t\in K^{\times}\Rightarrow n_\r(1)\ \hbox{acts as the identity on}\ L_\a.$$
This establishes that $\tilde\ll_{\b, \a}$ is determined, up to a
nonzero scalar multiple, by $\a$ and $\b$ and does not otherwise
depend on our choice of $w$ or $n$. Since also $n\dmg_\a=\dmg_\b, \
ne_\a=\vep e_\b$ for some $\vep\in K^{\times}$ (actually $\vep=\pm
1$).

\begin{defi}  $\ll_{\b, \a}: L_\a\to L_\b$ is the $K$-linear map
$\vep^{-1}\tilde\ll_{\b, \a}.$
\end{defi}

Clearly, $\ll_{\b, \a}$ depends only on $\a, \b$ and the choice of
Chevally basis. We see that
$$\ll_{\a, \b}=\ll_{\b, \a}^{-1}, \eqno(3.5)$$
$$\ll_{\a, \b}\ll_{\b, \r}=\ll_{\a, \r}, \eqno(3.6)$$
$$\ll_{\a, \a}=1. \eqno(3.7)$$

With $\a$ fixed as above, define $R=L_\a$ (as a $K$-space).
Eventually, $R$ will have a life of its own, independent of $\a$.
For this reason, given $r\in R$, we shall write it as $e_\a(r)$ when
we think of it as being in $L_\a$. We identify $K$ as a subspace of
$L_\a$ by $ae_\a=e_\a(a)$ for some $a\in K$.

Now if $\b\in \Delta $ is arbitrary, then the $K$-linear map
$\ll_{\b, \a}:L_\a\to L_\b$ is defined and we define
$$e_\b(r)=\ll_{\b, \a}e_\a(r).\eqno(3.8)$$
 From (3.7), we see that this definition is consistent when $\b=\a$.
From (3.5) and (3.6), we see that (3.8) holds for all pairs of roots
 $\a, \b\in\Delta.$

We note that for all $a, b\in K,\, r, s\in R$, we have
$e_\b(ar+bs)=ae_\b(r)+be_\b(s)$. This allows us to write $x_\b(r),\,
r\in R$, for any $x_\b\in \dmg_\b$ unambiguously. Clearly,
$$x_\b(r)=x_\b(s)\Longleftrightarrow r=s.$$

\begin{lemm}
For all $n\in N,\, \b\in\Delta,\, r\in R$, we have
$$(\Ad n)(e_\b(r))=((\Ad n)e_\b)(r).$$
\end{lemm}
{\bf Proof.} Let $nT\leftrightarrow w\in W$ and $\r=w\b$, so both
sides of the asserted equality lie in $L_\r$. Since $\ll_{\r, \b}$
equals, up to a scalar, $\Ad(n)|_{L_\b}$, $b\ll_{\r, \a}=\Ad
n(\ll_{\b, \a})$, for some $b\in K^{\times}$. Then
$$be_\r=b\ll_{\r, \a}e_\a=(\Ad n)\ll_{\b, \a}e_\a=(\Ad n)e_\b$$ and hence
$$(\Ad n)e_\b(r)=(\Ad n)\ll_{\b, \a}e_\a(r)=b\ll_{\r, \a}e_\a(r)=be_\r(r)
=((\Ad n)e_\b)(r).\quad \rule[-.23ex]{1.0ex}{2.0ex}$$

Now we set $$[e_\b(r), e_\r(s)]=[e_\b, e_\r](m_{(\b, \r)}(r,
s)).\eqno(3.9)$$

\begin{lemm}
If $(\b', \r')\in [(\b, \r)]$, then $m_{(\b', \r')}=m_{(\b, \r)}$.
\end{lemm}
{\bf Proof.}  Let $m'=m_{(\b', \r')}, m=m_{(\b, \r)}$. Choose $w\in
W$ with $\b'=w\b$ and $\r'=w\r$. Let $a, b, c, \vep, \vep'\in
K^{\times}$ be chosen such that
$$a\Ad(n)e_\b=e_{\b'}, \  b\Ad(n)e_\r=e_{\r'},\  c\Ad(n)e_{\b+\r}=e_{\b'+\r'},$$
$$[e_\b, e_\r]=\vep e_{\b+\r}, \  [e_{\b'}, e_{\r'}]=\vep' e_{\b'+\r'}.$$
Since $\Ad(n)$ is an automorphism, $\vep'=abc^{-1}\vep$. Now for
$r, s\in R$,
\begin{eqnarray*}
&&[e_{\b'}, e_{\r'}]m'(r, s)=[e_{\b'}(r), e_{\r'}(s)]\\
&=&[\ll_{\b', \b}e_{\b}(r), \ll_{\r', \r}e_{\r}(s)]\\
&=&[a\Ad(n)e_\b(r), b\Ad(n)e_\r(s)]=ab\Ad(n)[e_\b(r), e_\r(s)]\\
&=&ab\Ad(n)[e_\b, e_\r](m(r, s))\\
&=&ab\vep\Ad(n)e_{\b+\r}(m(r, s))=ab\vep c^{-1} e_{\b'+\r'}(m(r, s))\\
&=&\vep'e_{\b'+\r'}(m(r, s))=[e_{\b'}, e_{\r'}](m(r, s)).
\end{eqnarray*}
Thus  $m_{(\b', \r')}=m_{(\b, \r)}$  and $m=m'$.
\rule[-.23ex]{1.0ex}{2.0ex}

 From the above Lemma, we can define two multiplications on $R$:

(1) For a positive $A_2$-pair $(\b, \r)$ (see the Remark after Lemma
3.3), we define $\dashv: R\times R\to R$ given by
$$[e_\b(r), e_\r(s)]=[e_\b, e_\r](r\dashv s).$$

(2) For a negative $A_2$-pair $(\b, \r)$, we define $\vdash: R\times
R\to R$ given by
$$[e_\b(r), e_\r(s)]=[e_\b, e_\r](s\vdash r).$$

Clearly from Lemmas 3.5 and 3.6, we know that the above definition
is well-defined on $D$. Moreover, if $\Delta$ is  of type $D_l\;
(l\ge4)$ or $E_l\; (l=6, 7, 8)$, then $r\dashv s=s\dashv r$, i.e.,
$R$ is a commutative dialgebra since $m_{\b, \r}(r, s)=m_{\r, \b}(r,
s)$.

Let $(\a, \b)$ be a $A_2$-pair. Then $r_\b\a=\a+\b$. Suppose that
$e_\a, e_\b$ have already been chosen as our Chevalley basis. Let
$n_\b=\exp\ad e_\b\exp\ad (-e_{-\b})\exp\ad e_\b$. Then by
definition, we get that
\begin{eqnarray*}
n_\b e_\a(r)=n_\b(e_\a(r))&=&\exp\ad e_\b\exp\ad (-e_{-\b})\exp\ad e_\b(e_\a(r))\\
&=&\exp\ad e_\b\exp\ad (-e_{-\b})(e_\a(r)-[e_\a(r), e_\b])\\
&=&\exp\ad e_\b(e_\a(r)-[e_\a(r), e_\b]-[e_\a(r)+[e_\a(r), e_\b], e_{-\b}])\\
&=&-[e_\a(r), e_\b].
\end{eqnarray*}
Thus $$(n_\b e_\a)(r)=-[e_\a(r), e_\b].$$ In particular, $$n_\b
e_\a=-[e_\a, e_\b].$$

The element $e_\a\in L_\a$ defines an element of $R$ that has been
identified with $1\in K: e_\a(1)=1e_\a$. For any $\b\in \Delta$ we
have $e_\b(1)=\ll_{\b, \a}(e_\a(1))=e_\b$, so every basis element
$e_\b$ defines the same element of $R$.

Note that $1$ is bar-unit of $R$. For a positive $A_2$-pair $(\a,
\b)$, we have $[e_\a(r), e_\b(1)]=[e_\a, e_\b](r\dashv 1)$. But
\begin{eqnarray*}
[e_\a(r), e_\b(1)]&=&[e_\a(r), e_\b]=-n_\b(e_\a(r))\\
&=&-(n_\b e_\a)(r)=[e_\a, e_\b](r).
\end{eqnarray*}
Thus $r\dashv 1=r$. Similarly, for a negative $A_2$-pair $(\a, \b)$,
we can prove that $1\vdash r=r.$ So $R$ is a unital dialgebra.

Next, we investigate the associativity of the multiplication on $R$.

Now we assume that $l\ge 3$. Let $(\b, \r)$ and $(\r, \de)$ be two
positive $A_2$-pairs with $(\b, \de)=0$ (such pairs exist if
$l\ge3$). By Lemma 3.3 (3), $m_{(\b, \r)}=m_{(\r, \de)}$. By the
Jacobi identity and $[L_\b, L_\de]=0$, we obtain that
$$[[e_\b(r), e_\r(s)], e_\de(t)]=[e_\b(r), [e_\r(s),  e_\de(t)]],$$
$$[e_\b(r), [e_\de(t), e_\r(s)]]=-[[e_\b(r), e_\r(s)],  e_\de(t)],$$
$$[[e_\r(s), e_\b(r)], e_\de(t)]=-[[e_\r(r), e_\de(t)],  e_\r(s)].$$
It follows that $(r\dashv s)\dashv t=r\dashv (s\dashv
t)=r\dashv(s\vdash t)$ and $(r\vdash s)\dashv t=r\vdash (s\dashv
t)$. Similarly, we can prove the left identities in the axiom(Ass).
The $l=2$ case will be handled in \S5.2.

From the above, we have proved the following theorem.
\begin{theo}
Let $L$ be a $\Delta$-graded Leibniz algebra over a field $K$ of
characteristic $0$, where $\Delta$ is a simply-laced finite
indecomposable root system of rank $l\ge 2$. Let $\dmg$ be the
associated split simple Lie subalgebra with root system $\Delta$.
For any root $\a\in\Delta$ and let $R=L_\a$ as a $K$-vector space.
Relative to a Chevally basis $\{e_\b\mid \b\in\Delta\}\cup\{H_i\mid
i=1, \cdots, l\}$, define the map $\ll_{\b, \a}:L_\a\to L_\b$ of
Definition 3.4 and the element $e_\b(r), r\in R$ of (3.8). Then $R$
is a unital $K$-dialgebra. Moreover, $R$ is associative if the rank
$l\ge 3$. If $\Delta$ is of type $D$ or $E$, then $R$ is
commutative.
\end{theo}

\section{Centrally isogenous of $\Delta$-graded Leibniz algebra}

Let $L$ and $L'$ be $\Delta$-graded Leibniz algebras over $K$ for
the same finite root system $\Delta$. Then their associated split
simple subalgebras $\dmg$ and $\dmg'$ are isomorphic. Simply denoted
them by $\dmg$, so $\dmg$ is a subalgebra both of $L$ and $L'$.

\begin{defi}
A homomorphism $\vp: L\to L'$ is $\Delta$-homomorphic if
$\vp|_{\dmg}=id_{\dmg}$.
\end{defi}

Let $\vp: L\to L'$ be a $\Delta$-homomorphism. From the definition,
it follows at once that $\vp(L_\a)\subset L_\a'$ for all $\a\in
\Delta\cup\{0\}$. Let $R$ and $R'$ be the $K$-dialgebras associated
to $L$ and $L'$ respectively defined by a certain choice of
Chevalley basis for $\dmg$. For each $\a\in \Delta$, we have
$R=L_\a\stackrel{\vp}{\to}L_\a'=R'$, so $\vp$ determines a map
$$\bar\vp_\a: R\to R'.$$
\begin{prop}
$(1)$ $\bar\vp_\a: R\to R'$ is independent of the choice of
$\a\in\Delta$ (so we can denote it by $\bar\vp$) and is a
homomorphism of dialgebras. Furthermore, $\bar\vp$ is injective
(resp. surjective, bijective) if $\vp$ is injective (resp.
surjective, bijective).

$(2)$ If $\bar\vp$ is an isomorphism, then the $\Delta$-homomorphism
$\vp$ is central and $L$ and $L'$ are centrally isogenous.
\end{prop}

{\bf Proof.} (1) Let $\{e_\b\}_{\b\in\De}\cup \{H_i\}_{i=1}^l$ be
the Chevalley basis of $\dmg$. We have $L_\a=\{e_\a(r)\mid r\in R\}$
and $L_\a'=\{e_\a(r)\mid r\in R'\}$. Furthermore, $ke_\a=e_\a(k)$
for all $k\in K$ and $e_\a(k_1r+k_2s)=k_1e_\a(r)+k_2e_\a(s)$ for all
$r, s\in R, R'$, $ k_1, k_2\in K$. By definition,
$\vp(e_\a(r))=e_\a(\bar\vp_\a(r))$ for all $r\in R$, and it follows
that $\bar\vp_\a$ is $K$-linear and maps $1\in R$ to $1\in R'$.

Now suppose that $(\a, \b)$ is a positive $A_2$-pair. Then
\begin{eqnarray*}
[e_\a, e_\b](\bar\vp_{\a+\b}(r\dashv s))&=&\vp([e_\a, e_\b](r\dashv s))
=\vp([e_\a(r), e_\b(s)])\\
&=&[e_\a(\bar\vp_\a(r)), e_\b(\bar\vp_\b(s))]=[e_\a,
e_\b](\bar\vp_\a(r)\dashv\bar\vp_\b(s)).
\end{eqnarray*}
Thus  $$\bar\vp_{\a+\b}(r\dashv
s)=\bar\vp_\a(r)\dashv\bar\vp_\b(s).\eqno(4.1)$$ Similarly, we have
$$\bar\vp_{\a+\b}(r\vdash
s)=\bar\vp_\a(r)\vdash\bar\vp_\b(s).\eqno(4.2)$$

With $s=1$, we get $\bar\vp_{\a+\b}(r)=\bar\vp_\a(r)$. From this, it
is easy to see that $\bar\vp$ is independent of $\a$. Thus, $(4.1)$
and (4.2) show the homomorphism property of $\bar\vp$. The remaining
parts of (1) are obvious.

(2) If $\bar\vp$ is an isomorphism, then  $\bar\vp_\a$ is an
isomorphism for each $\a\in\De$, so $\ker\vp\subset L_0$. Since
$[L_\a, \ker\vp]\subset L_\a\cap\ker\vp=\{0\}$, we see that
$\ker\vp$ lies in the center of $L$ from Definition 3.1. Since $L$
and $L'$ are perfect and $L'\cong L/Z$ for some central ideal $Z$,
they have the same universal central extension.
\rule[-.23ex]{1.0ex}{2.0ex}

\begin{prop}
Let $L$ be a Leibniz algebra graded by $\Delta$ and $(\frak u, \vp)$
the universal central extension of $L$. Then $\frak u$ is graded by
$\Delta$ and has the same associated dialgebra as that of $L$.
Furthermore, $\vp$ is a $\Delta$-homomorphism and $\vp: \frak
u_\a\to L_\a$ is a homomorphism for all $\a\in\Delta$. In
particular, $\ker\vp\subset \frak u_0$.
\end{prop}
{\bf Proof.} It is well known that $\dmg$ is centrally closed in the
category of Leibniz algebras (\cite{Gao2}). Thus the central
extension $\vp: \vp^{-1}(\dmg)\to \dmg$ splits and we may view
$\dmg$ as a subalgebra of $\frak u$. In particular, $H$ is a
subalgebra of $\frak u$. We define
\begin{eqnarray*}
&&\tilde\frak u_\a:=\vp^{-1}(L_\a), \quad \a\in \De\cup\{0\},\\
&&\frak u_\a:=\cases{[\tilde\frak u_\a, H], \quad \a\in\De,\cr
\tilde\frak u_0, \quad \a=0.}
\end{eqnarray*}
For all $h_1, h_2\in H, x\in \tilde\frak u_\a$, we have
$$\ad h_1([x, h_2])=-[[x, h_2], h_1]=-[[x, h_1], h_2]=[\a(h_1)x+z, h_2]=\a(h_1)[x, h_2],$$
where $z\in\ker\vp$. This proves that
$$\frak u_\a \ \hbox{is an  $\a$-weight space for}\  \ad_{\frak u}H,\, \a\in\De.\eqno(4.3)$$
 It follows that for $\a\in\De$, $\frak u_\a\cap \ker\vp=\{0\}$ and hence
$$\vp|_{\frak u_\a}:\frak u_\a\to L_\a\eqno(4.4)$$ is an isomorphism of vector spaces.

Let $x\in\frak u$ and write $x=\sum_{\a\in\De\cup\{0\}}\tilde x_\a$,
where $\tilde x_\a\in\tilde\frak u_\a$. Fix $h\in H$ so that for all
$\a\in\De, \a(h)\ne 0$. Then
$$\tilde x_\a-\a(h)^{-1}[\tilde x_\a, x]\in\ker\vp\subset \tilde \frak u_0=\frak u_0.\eqno(4.5)$$
Thus we may rewrite $x$ as $\sum_{\a\in\De\cup\{0\}} x_\a$, where
$x_\a\in\frak u_\a$ and $x_0\in \frak u_0$. It follows that
$$\frak u=\frak u_0+\sum_{\a\in\De}\frak u_\a.$$

Now
$$\frak u_0=\tilde\frak u_0=\vp^{-1}(\sum_{\a\in\De}[L_\a, L_{-\a}])
=\sum_{\a\in\De}[\tilde \frak u_\a, \tilde \frak u_{-\a}]+\ker\vp.$$
Since $\tilde \frak u_\a=\frak u_\a+\ker\vp$ from (4.3), we have
$$\frak u_0=\sum_{\a\in\De}[\frak u_\a, \frak u_{-\a}]+\ker\vp.$$
This proves that $\frak u_0$ is a $0$-eigenspace for $H$. Now we see
that $\frak u_\a$ is exactly the $\a$-eigenspace for $H$ for all
$\a\in\De\cup\{0\}$ and $[\frak u_\a, \frak u_\b]\subset \frak
u_{\a+\b}$, whenever $\a+\b\in\De\cup\{0\}$. Also since $\frak
u=[\frak u, \frak u]$, we see that
$$\frak u_0=\sum_{\a\in\De}[\frak u_\a, \frak u_{-\a}]+[\frak u_0, \frak u_0].$$
But
$$[\frak u_0, \frak u_0] =\sum_{\a, \b\in\De}[[\frak u_\a, \frak u_{-\a}],
[\frak u_\b, \frak u_{-\b}]]\subset \sum_{\r\in\De} [\frak u_\r,
\frak u_{-\r}],$$
 so $$\frak u_0=\sum_{\a\in\De}[\frak u_\a, \frak
u_{-\a}].$$ This proves that $\frak u$ is graded by $\De$. From
(4.4) and the construction of $\vp$, we see that $\vp$ is a
$\De$-homomorphism and $\bar\vp$ is an isomorphism. Thus the
dialgebra associated to $\frak u$ is the same as that associated to
$L$. \rule[-.23ex]{1.0ex}{2.0ex}

\begin{prop}
Let $L$ be a Leibniz algebra graded by $\Delta$, $Z$ the center of
$L$ and $Z'$ any subspapce of $Z$. Then

$(1)$  $Z\subset L_0$.

$(2)$  $L/Z'$ has a unique structure as a Leibniz algebra graded
by $\Delta$ that makes the natural map $\pi: L\to L/Z'$ a
$\De$-homomorphism.

$(3)$  For all $\a\in\De$, $L_\a\stackrel{\pi}{\cong}(L/Z')_\a$ as
vector spaces.

$(4)$  The dialgebras associated to $L$ and $L/Z'$ are isomorphic
by the map $\bar\pi$ induced by $\pi$ as in Proposition 4.2.
\end{prop}
{\bf Proof.} Since for all $h\in H$, $\ad h|_{L_\a}$ is a scalar by
$\langle \a, h \rangle$, we see that $Z\subset L_0$. The remaining
results are now obvious. \rule[-.23ex]{1.0ex}{2.0ex}

\begin{prop}
Let $L$ and $L'$ be centrally isogenous Leibniz algebras and suppose
that $L$ is graded by a finite root system $\Delta$. Then $L'$ is
also graded by $\Delta$ and in such a way that the associated
dialgebras are isomorphic.
\end{prop}
{\bf Proof.} Let $\frak u$ be a universal central extension of $L$.
By Proposition 4.3, $\frak u$ is graded by $\De$. By assumption,
$L'\cong \frak u/Z_1$ for some subspace of the center of $\frak u$.
By Proposition 4.4, $L'$ is graded by $\De$. The associated
dialgebras are isomorphic by Propositions 4.3 and 4.4.
\rule[-.23ex]{1.0ex}{2.0ex}

In short, all Leibniz algebras in a given isogeny class are
$\De$-graded if one of them is, and all have isomorphic root spaces
for all $\a\in\De$. They differ only by central elements in the
$0$-root space.

\section{Proof of Recognition Theorem}

In this section we shall complete the proof of the Recognition
Theorem. We use the results and notation of $\S 3$ and $\S4$ freely.

\subsection{The Recognition Theorem for types $D$ and $E$}

Now we assume that $L$ is a Leibniz algerba graded by $\De$, where
$\De$ is a finite root system of type $D_l$, $l\ge 4$, or $E_6, E_7,
E_8$.

\begin{lemm} Let $\a, \b\in\De, r, s \in R $ and set $\rho=:\ad([e_\a(r), e_{-\a}(s)])$. Then
$$\rho(e_\b(t))=\langle \b, \a^{\vee}\rangle e_\b((t\dashv r)\dashv s).$$
\end{lemm}
{\bf Proof.} Case 1. $(\a|\b)=-1$:  Then $\b-\a\not\in\De$ and we
have
\begin{eqnarray*}
\rho(e_\b(t))&=&-[e_\b(t), [e_\a(r), e_{-\a}(s)]]\\
             &=&-[[e_\b(t), e_\a(r)], e_{-\a}(s)]\\
             &=&-[[e_\b, e_\a](t\dashv r), e_{-\a}(s)]\\
             &=&-[[e_\b, e_\a], e_{-\a}]((t\dashv r)\dashv s)\\
             &=&-[[e_\b, [e_\a, e_{-\a}]]((t\dashv r)\dashv s)\\
             &=&\langle \b, \a^{\vee}\rangle e_\b((t\dashv r)\dashv
             s).
\end{eqnarray*}

Case 2. $(\a|\b)=1$: Then $(-\a|\b)=-1$ and we may use Case 1.

Case 3. $(\a|\b)=0$: Then neither $\b+\a$ nor $\b-\a$ is a root
and both sides of our equation are 0.

Case 4. $\a=\b$:  It is easy to find a pair of roots $\r,
\vep\in\De$ with $\r+\vep=\a$, $\r, \vep\not\in\{\a, -\a\}$. Then
$e_\b(t)=e_\a(t)=[e_\r, e_\vep](t')$ for some $t'\in R$ which is
some $K$-multiple of $t$. Applying $\rho$ and using the previous
cases,
\begin{eqnarray*}
\rho(e_\b(t))&=&\rho([e_\r, e_\vep](t'))=\rho([e_\r(t'), e_\vep(1)])\\
             &=&\langle \r, \a^{\vee}\rangle [e_\r((t'\dashv r)\dashv s), e_\vep(1)]+\langle \vep, \a^{\vee}\rangle [e_\r(t'), e_\vep(r\dashv s)]\\
             &=&\langle \r+\vep, \a^{\vee}\rangle [e_\r, e_\vep](t'\dashv(r\dashv s))\\
             &=&\langle \a, \a^{\vee}\rangle e_\b(t\dashv (r\dashv s))\\
             &=&\langle \a, \a^{\vee}\rangle e_\b((t\dashv r)\dashv
             s).
\end{eqnarray*}

Case 5. $-\a=\b$: This is similar to Case 4.
\rule[-.23ex]{1.0ex}{2.0ex}

We wish to define a homomorphism
$$\phi:L\to \frak t(R, \De)=\dmg\ot R\eqno(5.1)$$
by defining
$$\phi(e_\a(r))=e_\a\ot r, \quad \hbox{for all}\ \a\in\De, r\in R.\eqno(5.2)$$

Since $L$ is generated by the subspace $L_\a, \a\in\De$, it is clear
that (5.2) uniquely defines $\phi$.

Suppose that $\sum_{\a\in\De}\sum_{i=1}^{n_\a}[e_\a(r(\a, i)),
e_{-\a}(s(\a, i))]=0$ for some $r(\a, i), s(\a, i)\in R$. Then using
Lemma 5.1, we have for all $\b\in\De$
$$\sum_{\a\in\De}\sum_{i=1}^{n_\a}\langle \b, \a^{\vee}\rangle e_\b(r(\a, i)\dashv
s(\a, i))=0$$
and hence
$$\sum_{\a\in\De}\langle \b, \a^{\vee}\rangle\sum_{i=1}^{n_\a} e_\b(r(\a, i)
\dashv s(\a, i))=0.$$ Since $e_\b(r)=e_\b(s)\Leftrightarrow r=s$. It
follows that in $H\ot R$,
$\sum_{\a\in\De}\sum_{i=1}^{n_\a}\a^{\vee}\ot(r(\a, i)\dashv s(\a,
i))=0$ (it is killed by all the functionals $1\ot \b$). Thus we may
extend $\phi$ to $L_0$ by
$$\phi: \sum_{\a\in\De}\sum_{i=1}^{n_\a}[e_\a(r(\a, i)), e_{-\a}(s(\a, i))]\mapsto \sum_{\a\in\De}\sum_{i=1}^{n_\a}\a^{\vee}\ot (r(\a, i)\dashv s(\a, i)) .$$

We check that $\phi$ is a homomorphism: If $(\a, \b)$ is an
$A_2$-pair, then
\begin{eqnarray*}
\phi([e_\a(r), e_\b(s)])&=&\phi([e_\a, e_\b](r\dashv s))\\
                        &=&[e_\a, e_\b]\ot(r\dashv s)\\
                        &=&[\phi(e_\a(r)), \phi(e_\b(s))].
\end{eqnarray*}

If $(\a, \b)\ge 0$, then $\phi([e_\a(r),
e_\b(s)])=0=[\phi(e_\a(r)), \phi(e_\b(s))]$.

If $\b=-\a$ then  $\phi([e_\a(r),
e_{-\a}(s)])=\a^{\vee}\ot(r\dashv s)=[e_\a\ot r, e_\b\ot s]$. Now
let $h=[e_\a(r), e_{-\a}(s)]$. We have
\begin{eqnarray*}
\phi([e_\b(t), h])&=&-\phi(\langle \b, \a^{\vee}\rangle e_\b((t\dashv r)\dashv s))\\
                  &=&-\langle \b, \a^{\vee}\rangle e_\b\ot ((t\dashv r)\dashv s)\\
                  &=&[e_\b\ot t, \a^{\vee}\ot (r\dashv s)]\\
                  &=&[\phi(e_\b(t)), \phi(h)],
\end{eqnarray*} and
\begin{eqnarray*}
&&\phi([[e_\b(t), e_{-\b}(u)], h])\\
&=&-\phi(\langle \b, \a^{\vee}\rangle [e_{\b}(t), e_{-\b}((u\dashv r)\dashv s)]+\langle \b, \a^{\vee}\rangle [e_\b((t\dashv r)\dashv s), e_{-\b}(u)])\\
                  &=&-\langle \b, \a^{\vee}\rangle \Big(\b^{\vee}\ot (t\dashv ((u\dashv r)\dashv s))-\b^{\vee}\ot (((t\dashv r)\dashv s)\dashv u)\Big)\\
                  &=&0=[\phi([e_\b(t), e_{-\b}(u)]), \phi(h)]
\end{eqnarray*}
since $ ((t\dashv r)\dashv s)\dashv u=u\vdash((t\dashv r)\dashv
s)=(u\vdash(t\dashv r))\dashv s=((u\vdash t)\dashv r)\dashv
s=((t\dashv u)\dashv r)\dashv s=t\dashv ((u\dashv r)\dashv s)$,
where we have used the commutativity and associativity of $R$.

\begin{prop}
The homomorphism $\phi$ define by (5.1) and (5.2) is a surjective
$\De$-homomorphism and $\ker\vp$ is contained in the center of
$L$.
\end{prop}
{\bf Proof.} Since $L_\a=e_\a(R)\stackrel{\phi}{\to}e_\a\ot R$, it
is clear that $\phi$ is bijective on the root spaces of $L_\a,
\a\in\De$ and hence $\phi$ is surjective and $\ker\vp\subset L_0$.
Thus $[\ker\phi, L_\a]\subset L_\a\cap (\ker\phi)=(0)$ and hence
$\ker\phi$ is central by (2) of Definition 3.1.
\rule[-.23ex]{1.0ex}{2.0ex}

We have proved that $L$ is a central extension of $\dmg\ot R$, thus
proving the third part of the Recognition Theorem 1.1.  The
universal central extension of $\dmg\ot R$ is given in \cite{LL} for
a unital commutative associative dialgebra $R$.

\subsection{The Recognition Theorem for type $A$}

Now we assume that $L$ is a Leibniz algerba graded by $\De$, where
$\De$ is a finite root system of type $A_l$, $l\ge 2$.

Let $n=l+1$ and $\vep_1, \cdots, \vep_n$ be an orthonormal basis
for $\mathbb R^n$. Identify $\Delta$ with $\{\vep_i-\vep_j\mid
i\ne j\}$ and define a base $\{\a_1, \cdots, \a_l\}$ of $\De$ by
$\a_i=\vep_i-\vep_{i+1}$. The positive class of $A_2$-pair in
$\De$ will be taken as $[(\a_1, \a_2)]$.

For $\a=\vep_i-\vep_j\in\De$ we let $L_{ij}=L_\a$. The simple Lie
algebra $\dmg$ over $K$ of type $A_l$ may be identified with
$\fsl(n, K)$. We choose as our Chevalley basis the matrix units
$E_{ij}, i\ne j$ and the elements $h_i=[E_{i, i+1}, E_{i+1, i}],
i=1, \cdots, l$. To bring the notation in line with $\S2.4$, we
write $e_{ij}$ for $E_{ij}$. Let $R$ be the $K$-dialgebra derived
from $L$ with this choice of positive $A_2$-pairs and the given
Chevalley basis. Then we have
$$[(\a_1, \a_2)]=\{\vep_i-\vep_j, \vep_j-\vep_k\mid i, j, k\ \hbox{distinct}\}$$ and
$$[(\a_2, \a_1)]=\{\vep_i-\vep_j, \vep_k-\vep_i\mid i, j, k\ \hbox{distinct}\}.$$

Thus by results in $\S4$, we have
\begin{eqnarray}
&(i)&  L_{ij}=K\{e_{ij}(r)\mid r\in R\}, \,i\ne j;\\
&(ii)& e_{ij}(k_1r+k_2s)=k_1e_{ij}(r)+k_2e_{ij}(s) ;\\
&(iii)& [e_{ij}(r), e_{kl}(s)]=0\  \hbox{ if }\ i\ne l\ \hbox{ and } \  j\ne k;\\
&(iv)& [e_{ij}(r), e_{kl}(s)]=e_{il}(r\dashv s)\  \hbox{ if } \ i\ne l\ \hbox{and } \  j= k;\\
&(v)& [e_{ij}(r), e_{kl}(s)]=-e_{kj}(s\vdash r)\  \hbox{ if } \ i=l\
\hbox{and} \  j\ne k,
\end{eqnarray}
for all  $r, s\in R, \ k_1, k_2\in K$.

Now whether $R$ is associative or not, $L$ is homomorphic to the
image of $\stl(n, R)$, under the map
$$\phi: v_{ij}(r)\to e_{ij}(r), \quad r\in R,\, i\ne j.\eqno(5.6)$$

In Section 3, we show that $R$ is associative if $l\ge3$. Similar to
the proof of Proposition 3.1 in \cite{L}, we can show that $R$ is
alternative if $l=2$.

\begin{prop}
The homomorphism of $(5.6)$ is central.
\end{prop}
{\bf Proof.} It is clear that $\stl(n, R)$ is graded by $\De$ and
its associated dialgebra is $R$. The homomorphism (5.6) sends
$v_{ij}(R)\to e_{ij}(R)=L_{ij}$ isomorphically and hence the
$\ker(\phi) \subset \stl(n, R)_0$.  The same argument as used in
Proposition 5.2 gives that $\phi$ is central.
\rule[-.23ex]{1.0ex}{2.0ex}

This completes (1) and (2) of the Recognition Theorem 1.1.

\noindent{\bf Remark.} 1. The structures of Leibniz algebras graded
by finite root systems of other types are determined in \cite{LH2}
by using the methods in \cite{BZ} and \cite{ABG2}.

2. By Theorem 1.3 in \cite{MZ} and Theorem 1.1, we can easily obtain
the following theorem.
\begin{theo}
 A $Q(n)$-graded Leibniz superalgebra over $K$ is centrally isogenous to $\stl(n+1, D)$,
 where $D=D_{\bar0}+D_{\bar1}$ is an associative or an alternative (if $n=2$)
 unital dialgebra such that there exists $\nu\in D_{\bar1},
 \nu^2=1$. (The definition of $Q(n)$-graded Leibniz superalgebra immediately follows
 from Definition 1.3 in \cite{MZ} and
 Definition 3.1 in section 3, also see \cite{LH3}).
\end{theo}

\vskip45pt \centerline{\bf ACKNOWLEDGMENTS}

\vskip15pt    Part of the work was done during the first author's
visit in the Department of Mathematics at University of Bielefeld
and study in the Department of Mathematics at Shanghai Jiaotong
University for his Post-doctor study. He expresses his gratitude to
the `AsiaLink Project' for financial support. He is also indebted to
C.M. Ringel and C.P. Jiang for the hospitality and encouragement.
Authors give their thanks to Y. Gao for his useful comments. Liu is
supported by the NNSF (Grant 10671027, 10701019 and 10571119), the
ZJNSF(Grant Y607136), and Qianjiang Excellence Project of Zhejiang
Province (No. 2007R10031).  Hu is supported in part by the NNSF
(Grants 10431040, 10728102), the PCSIRT, the Priority Academic
Discipline from the MOE of China,  the Shanghai Priority Academic
Discipline from the SMEC. Authors are grateful to the referee for
correction in some errors and invaluable suggestions.

\vskip30pt
\def\refname{\cen{\normalsize\bf REFERENCES}}

\end{document}